\documentclass{amsart}
\usepackage[all]{xy}
\usepackage{amsfonts,amscd,amssymb,amsmath,mathrsfs}
\theoremstyle{plain}
\newtheorem{thm}{Theorem}[section]

\newtheorem{prop}[thm]{Proposition}
\newtheorem{cor}[thm]{Corollary}
\theoremstyle{definition}
\newtheorem{defn}{Definition}[section]
\newtheorem{ex}[defn]{Example}
\theoremstyle{remark}
\newtheorem{rmk}{Remark}[section]

\numberwithin{equation}{section}

\begin{document}

\newcommand{\Soc}{\operatorname{Soc}}
\newcommand{\modP}{\mod_{P}\Lambda}
\newcommand{\modl}{\mod \Lambda}
\newcommand{\module}{\operatorname{mod}}
\newcommand{\p}{\operatorname{p}}
\newcommand{\inj}{\operatorname{i}}
\newcommand{\Module}{\operatorname{Mod}}
\newcommand{\Cok}{\operatorname{Coker}}
\newcommand{\Hom}{\operatorname{Hom}}
\newcommand{\h}{\operatorname{h}}
\newcommand{\e}{\operatorname{e}}
\newcommand{\res}{\operatorname{res}}
\newcommand{\Tr}{\operatorname{Tr}}
\newcommand{\TrD}{\operatorname{TrD}}
\newcommand{\rad}{\operatorname {{\bold r}}}
\newcommand{\La}{\operatorname{\Lambda}}
\newcommand\End{\operatorname{End}}
\newcommand\Ext{\operatorname{Ext^1_\Lambda}}
\newcommand\Ex{\operatorname{Ext}}
\newcommand\ann{\operatorname{ann}}
\newcommand\coend{\operatorname{coend}}
\newcommand\Img{\operatorname{Im}}
\newcommand\D{\operatorname{D}}
\newcommand\DTr{\operatorname{DTr}}
\newcommand\Ker{\operatorname{Ker}}
\newcommand\Coker{\operatorname{Coker}}
\newcommand{\tr}{\operatorname {t}}
\newcommand{\ct}{\operatorname{ct}}
\newcommand{\rej}{\operatorname{rej}}
\newcommand{\gd}{\operatorname{{gl.dim}}}
\newcommand{\radic}{\operatorname{rad}}
\newcommand{\add}{\operatorname{add}}
\newcommand{\Ind}{\operatorname{Ind}}
\newcommand{\Sub}{\operatorname{Sub}}
\newcommand{\Fac}{\operatorname{Fac}}
\newcommand{\lap}{\operatorname{l}}
\newcommand{\rap}{\operatorname{r}}
\newcommand{\locn}{\operatorname{lnRep}}
\newcommand{\rep}{\operatorname{Rep}}

\newcommand{\gl}{(\Gamma,\Lambda)}
\newcommand{\poset}{(\Gamma_{0}¥,\Lambda)}
\newcommand{\Stilde}{\tilde S}
\newcommand{\Ttilde}{\tilde T}
\newcommand{\ad}{(+){\rm -admissible}}
\newcommand{\pa}{{\tilde{\mathcal P}}(A)}
\newcommand{\modA}{{\rm mod}A}
\newcommand{\hull}{H_{\Lambda}}
\newcommand{\ds}[1]{\ensuremath{#1}}
\newcommand{\quiver}[2]{\ensuremath{\left( #1, #2 \right) }}
\newcommand{\quiverA}[1]{\ensuremath{\left( A_n, #1 \right) }}
\newcommand{\set}[1]{\ensuremath{\left\{#1\right\}}}

\bibliographystyle{plain}
\title[Admissible sequences and the Weyl group]{Admissible sequences, preprojective modules, and reduced words in the Weyl group of a quiver}
\author{Mark Kleiner and Allen Pelley}
\address{Department of Mathematics, Syracuse University, Syracuse, New 
York 13244-1150}
\email{mkleiner@syr.edu}
\email{anpelley@syr.edu}
\thanks{The authors are supported by the NSA grant H98230-06-1-0043.}

\begin{abstract} This paper studies connections between the preprojective modules over the 
path algebra of a finite connected quiver without oriented cycles, the
(+)-admissible sequences of vertices, and the Weyl group.  For each preprojective   module, there exists a unique up to a certain equivalence shortest (+)-admissible sequence annihilating the module.  A (+)-admissible sequence is the shortest sequence annihilating some preprojective module if and only if the product of  simple reflections associated to the vertices of the sequence  is a reduced word in the Weyl group.  These statements have the following application that strengthens known results of Howlett and Fomin-Zelevinsky.  For any fixed Coxeter element of the Weyl group associated to an indecomposable symmetric generalized Cartan matrix, the group is infinite if and only if the powers of the element are reduced words. 
\end{abstract}
\maketitle

\section*{Introduction}\label{Section:intro}

A preprojective module over the path algebra of a finite connected quiver without oriented cycles (or a (+)-irregular representation of the quiver)  was defined by Bernstein, Gelfand, and Ponomarev ~\cite{bgp}  as a module that can be annihilated (reduced to zero) by a finite sequence of operations, where each operation consists in choosing a sink (vertex at which no arrow starts), reversing the direction of each arrow ending at the sink, and taking the image of the module under  a suitable functor (positive reflection functor) into the category of modules over the path algebra of  the new quiver.  A sequence of vertices of the original quiver for which the indicated sequence of operations is possible is called a (+)-admissible sequence.  Let $\mathfrak S$ be the set of (+)-admissible sequences.  If $M$ is an indecomposable preprojective module, let $S_M$ be a shortest sequence in $\mathfrak S$ that annihilates $M$.  The sequence $S_M$ is unique (up to a certain equivalence $\sim$), and it can be constructed by a  simple combinatorial procedure if the location of $M$ in the preprojective component (see ~\cite{ars}) of the Auslander-Reiten quiver is known.  Conversely, $M$ is uniquely (up to isomorphism) determined by $S_M$.  These and other properties of $S_M$ were studied in ~\cite{kt} in order to get new insights into the structure of the preprojective component.  

In this paper we study connections between the category  $\tilde{\mathscr P}$ of preprojective modules, the set $\mathfrak S$, and  the Weyl group $\mathcal W$ ~\cite{bgp} of the underlying (nonoriented) graph of the quiver.  There were several indications in favor of undertaking such a study.  Although the authors of ~\cite{kt} did  not mention it explicitly, they had an interest in studying the elements of $\mathcal W$ associated to the sequences $S_M$.  Wolfgang Rump noted that ~\cite[Theorem 3.1]{kt}, which says that $S_M$ is unique up to equivalence and determines $M$, might admit an interesting formulation in terms of the group $\mathcal W$.  Then Andrei Zelevinsky suggested a procedure that should produce all indecomposable modules in $\tilde{\mathscr P}$ (up to isomorphism) from certain reduced words in $\mathcal W$.  One of our results is that the procedure works (Corollary \ref{zel}).

In order to carry out our study of $\tilde{\mathscr P}$, $\mathfrak S$, and  $\mathcal W$, we extend a result from  ~\cite{kt} by showing that for each  $M\in\tilde{\mathscr P}$ (not necessarily indecomposable), a shortest   sequence $S_M\in \mathfrak S$ annihilating $M$ is unique up to equivalence (Theorem \ref{shrtst}), and consider, for each $S\in\mathfrak S$,  the element $w(S)\in \mathcal W$ that is the composition of  simple reflections associated to the vertices of $S$.  
Using  properties of  $\mathfrak S$  and $\mathcal W$, we get information about $\tilde{\mathscr P}$.  For instance, if $M,N\in\tilde{\mathscr P}$ are indecomposable, then $M\cong N$ if and only if $w(S_M)=w(S_N)$ (Theorem \ref{reducedword}).  For all $S\in\mathfrak S$, there exists an $M\in\tilde{\mathscr P}$ satisfying $S\sim S_M$ if and only if the word $w(S)\in \mathcal W$ is reduced, and a simple procedure determines whether $w(S)$ is reduced (Theorems \ref{princseqreducedwrd} and \ref{seqreducedwrd}).  If $\Gamma$ is not a Dynkin diagram of the type A, D, or E, then for all $S\in\mathfrak S$,   the word $w(S)$ is reduced and there exists an $M\in\tilde{\mathscr P}$ satisfying $S\sim S_M$ (Corollary \ref{reduced}). We also give an elementary proof of the following well known statement (see \cite{ars}): if $M\in\tilde{\mathscr P}$  is indecomposable and if $N$ is an indecomposable module with the same dimension vector as that of $M$, then $M\cong N$ (Proposition \ref{reduced1}).  

Conversely, using properties of $\tilde{\mathscr P}$ and $\mathfrak S$, we obtain  information about $\mathcal W$.  Let $\mathcal W$ be a Coxeter group generated by  reflections $\sigma_1,\dots,\sigma_n$, and let  $c$ be any Coxeter element of  $\mathcal W$, i.e., $c=\sigma_{x_n}\dots\sigma_{x_1}$  where $x_1,\dots,x_n$ is any permutation of the  numbers $1,\dots,n$.  If  $A=(a_{ij})$ is an indecomposable  generalized $n\times n$ Cartan matrix and $\sigma_1,\dots,\sigma_n$ are the simple reflections, denote by ${\mathcal W}(A)$  the Weyl group ~\cite{kac1990}.  Zelevinsky brought to our attention the following two results.   Howlett proved that  $\mathcal W$ is  infinite if and only if $c$ has infinite order ~\cite[Theorem 4.1]{h}.  Fomin and Zelevinsky proved the  following. Let $A$ be symmetrizable and bipartite, i.e., the set $\{1,\dots,n\}$ is a disjoint union of nonempty subsets $I,J$ and, for $h\ne l$,  $a_{hl}=0$ if either $h,l\in I$ or $h,l\in J$.  If $c=\prod_{i\in I}\sigma_i\prod_{j\in J}\sigma_j$, then ${\mathcal W}(A)$ is infinite if and only if the powers of $c$ are reduced words in the $\sigma_h$'s ~\cite[Corollary 9.6]{CAIV}.  Inspired by the latter, we prove that if  $A$ is symmetric and $c$ is any Coxeter element, then ${\mathcal W}(A)$ is infinite if and only if the powers of $c$ are reduced words   (Theorem \ref{infWeylGrp}), which in our setting  strengthens the aforementioned results of Howlett and Fomin-Zelevinsky.  After we informed Robert Howlett of our result, he pointed out that Daan Krammer considered in  ~\cite{krammer94} elements all of whose powers are reduced words.

Our subsequent work will address the case of a symmetrizable matrix $A$.  For that we have to use representations of valued quivers studied by Dlab and Ringel ~\cite{dr} and an extension of the results of ~\cite{kt}  from quivers to valued quivers ~\cite{kt2}.

Another statement we prove is that if $M\in\tilde{\mathscr P}$, then  $w(S_M)^{-1}$  is a Coxeter-sortable element of $\mathcal W$ (see Definition \ref{c-sort1} and Proposition \ref{c-sort}).  Nathan Reading  introduced Coxeter-sortable elements for an arbitrary Coxeter group $\mathcal W$ ~\cite{rea} and  proved that if $\mathcal W$ is finite, then the set of Coxeter-sortable elements maps bijectively onto  the set of clusters ~\cite{FZ2003, mrz2003} and onto the set of noncrossing partitions ~\cite{McCam}.  It follows that to each preprojective module correspond a cluster and a noncrossing partition.  It would be interesting to investigate the clusters and noncrossing partitions thus obtained.

We now describe the content of the paper section by section. In 
Section 1 we recall the definitions and notation needed for the rest of the paper.  In Section 2 we develop combinatorial properties of $\mathfrak S$ needed in Sections 3 and 4.   The preorder $\preccurlyeq$  defined for $S,T\in\mathfrak S$ by setting $S\preccurlyeq T$ if $T\sim SU$  for some $U$ ~\cite{kt} induces a lattice structure on the set of equivalence classes of $\sim$.  The relation between the meet and the join in this lattice is similar to the relation between the greatest common divisor and the least common multiple in the set of integers (Theorem \ref{Theorem:meetjoin}).  The binary relations $\preccurlyeq$ and $\sim$ have the left cancellation property with respect to concatenation on $\mathfrak S$.  In Section 3, among other things, we obtain  two  crucial properties of principal (+)-admissible sequences (see Definition \ref{tight}).  If $S=x_1,\dots,x_s,\ s>1$, is a principal (+)-admissible sequence, then $T=x_2,\dots,x_{s}$ is a principal (+)-admissible sequence with respect to the orientation obtained by reversing the direction of each arrow ending at ${x_1}$ (Proposition \ref{princseq}).  If $M, N\in\tilde{\mathscr P}$ are  indecomposable and $[M],[N]$ are their isomorphism classes, then $S_M\preccurlyeq S_N$ if and only if there exists a path in the preprojective component starting at $[M]$ and ending at $[N]$ (Proposition \ref{fullsubposet}).  The latter statement illustrates the utility of (+)-admissible sequences by expressing a  complicated relation between indecomposable preprojective modules in terms of a simple relation between their shortest annihilating sequences.  Section 4 contains the main results and Section 5 deals with Coxeter-sortable elements.

One can obtain the results analogous to those of this paper by replacing preprojective modules ((+)-irregular representations) with preinjective modules (($-$)-irregular representations), and (+)-admissible sequences with ($-$)-admissible sequences ~\cite{bgp}.  We leave this to the reader.

The authors are grateful to Andrei Zelevinsky for many helpful suggestions and stimulating conversations, and to Nathan Reading for a brief introduction to Coxeter-sortable elements.

\section{Preliminaries}

We begin by recalling some facts, 
definitions, and notation, using freely ~\cite{ars,bgp}. 

A {\it graph} is a pair $\Delta=(\Delta_{0},\Delta_{1})$  where $\Delta_{0}$ 
is the set of vertices and $\Delta_{1}$ is the set of (possibly, 
multiple) edges of $\Delta$.  An orientation, $\Theta$,
on $\Delta$ consists of two functions $s:\Delta_{1}\to\Delta_{0}$ and
$e:\Delta_{1}\to\Delta_{0}$.  For an edge $a\in\Delta_{1}$, 
$s(a)$ and $e(a)$ are the 
vertices incident with $a$, and they are called the {\it starting point} and 
the
{\it endpoint} of $a$, respectively; one writes $a:s(a)\to e(a)$. The ordered pair 
$(\Delta,\Theta)$ is called a {\it quiver}, and $a$
is then called an arrow of $(\Delta,\Theta)$.  Given a sequence of arrows 
$a_{1},\dots,a_{t},\ t>0,$ satisfying $e(a_{i})=s(a_{i+1}),\ 
0<i<t,$ one forms a {\it path} $p=a_{t}\dots a_{1}$ of 
length $t$ in $(\Delta,\Theta)$.  By definition, $s(p)=s(a_{1}),\ e(p)=e(a_{t}),$ 
so one writes $p:s(p)\to e(p)$ and says that $p$ is a path from 
$s(p)$ to $e(p)$.  By definition, for all $x\in\Delta_{0}$ there is a 
unique path of length $0$ from $x$ to $x$. A path $p$ of length at least 1 is an {\it oriented cycle} if 
$s(p)=e(p)$. The set of vertices of any quiver without 
oriented 
cycles (no finiteness 
assumptions) acquires a structure of a 
partially ordered set (poset) by putting $x\le y$ if there is a 
path from $x$ to $y$.  If $(\Delta,\Theta)$ has no oriented cycles, we 
denote this poset by $(\Delta_0,\Theta)$.  For $x\in\Delta_{0}$, let $\sigma_{x}\Theta$ 
be the orientation on $\Delta$ obtained by reversing the
direction of each arrow incident with $x$ and preserving the 
directions of the remaining arrows.  There results a
new quiver $(\Delta,\sigma_{x}\Theta)$.  A vertex $x$ is a
{\it sink} (respectively, {\it source})  if no arrow starts (respectively, ends) at $x$. A sequence of vertices $S=x_{1},x_{2},\dots,x_{s}$ is called {\it (+)-admissible}
on $(\Delta,\Theta)$ 
if it either is  empty, or satisfies the following conditions: $x_{1}$ is a sink 
with respect to $\Theta$, $x_{2}$ is
a sink with respect to $\sigma_{x_{1}}\Theta$, and so on; we put 
$\Theta^{S}=\sigma_{x_{s}}\dots\sigma_{x_{1}}\Theta$.  If $T=u_1,\dots,u_p,v_1,\dots,v_q$ is (+)-admissible on  $(\Delta,\Theta)$, then $U=u_1,\dots,u_p$ is (+)-admissible on  $(\Delta,\Theta)$ and $V=v_1,\dots,v_q$ is (+)-admissible on  $(\Delta,\Theta^U)$; we write $T=UV$.

     Throughout the paper, $k$ is an arbitrary field, and 
$\Gamma=(\Gamma_{0}¥,\Gamma_{1})$ is a fixed finite 
connected graph without loops and with more than one vertex.
All orientations $\Lambda,\Theta$, etc., on $\Gamma$ are 
such that $\gl$, $(\Gamma,\Theta)$, etc., have no oriented cycles.  

A {\it representation} $(V,f)$ of a 
quiver $\gl$ over $k$ is a
set of finite-dimensional $k$-spaces $\{V(x)\,|\,x\in \Gamma_{0}\}$ together with
$k$-linear maps $f_{a}:V(x)\to V(y)$ for each arrow $a:x\to y$. 
We
denote by Rep$\gl$ the category of representations of $\gl$
over $k$, and by f.d.$\,k\gl$ the category of left 
modules of finite $k$-dimension over the (finite-dimensional) path algebra $k\gl$ (see ~\cite{ars}). In this 
paper all $k\gl$-modules belong to $\mathrm{f.d.}\,k\gl$. The 
categories Rep$\gl$ and $\mathrm{f.d.}\,k\gl$ are equivalent, and we view the
equivalence as identification. If $M\in\mathrm{f.d.}\,k\gl$ is identified with $(V,f)\in\,$Rep$\gl$, we define their {\it support} by $\mathrm{Supp}\,M=\mathrm{Supp}\,(V,f)=\{x\in\Gamma_0\,|\,V(x)\ne0\}$.  The {\it dimension vector} is $\dim M=\dim (V,f)=(\dim_k V(x))$ where $x\in \Gamma_0$.
 
 For each sink $x$ in $\gl$, the positive
reflection functor 
$F^{+}_{x}:\mathrm{f.d.}\,k\gl\to\mathrm{f.d.}\,k(\Gamma,\sigma_{x}\Lambda)$ 
is defined ~\cite[Definition 1.1, part 1)]{bgp}, and we recall the definition for the convenience of the reader.
If $M\in\mathrm{f.d.}\,k\gl$ is identified with $(V,f)\in\,$Rep$\gl$, let $(W,g)\in\,$Rep$(\Gamma,\sigma_{x}\Lambda)$ be identified with $F^{+}_{x}M$.  Then $W(z)=V(z)$ for all $z\in\Gamma_0\setminus\{x\}$, and $g_b=f_b$ for all those arrows $b$ of $(\Gamma,\sigma_{x}\Lambda)$ that do not start at $x$.  Let $a_i:y_i\to x,\ i=1,\dots,l$, be the arrows of $\gl$ ending at $x$, then the reversed arrows $a_i':x\to y_i,\ i=1,\dots,l$, are all the arrows of $(\Gamma,\sigma_{x}\Lambda)$ starting at $x$.  Consider the exact sequence
$$
0\to\Ker h\overset{j}\to\overset{l}{\underset{i=1}\oplus}V(y_i)\overset{h}\to V(x)
$$
of $k$-spaces, where the map $h$ is induced by the maps $f_{a_i}:V(y_i)\to V(x)$.  Then $W(x)=\Ker h$ and the maps $g_{a_i'}:W(x)\to W(y_i)=V(y_i)$ are  induced by $j$.
The functor $F^{+}_{x}$ is defined on morphisms in a natural way.

Replacing a sink with a
source  and a kernel with a cokernel, one gets a similar definition of a negative reflection functor $F_{x}^{-}$ ~\cite[Definition 1.1, part 2)]{bgp}.

We quote ~\cite[Theorem 1.1]{bgp}.
\begin{thm}\label{PosReflFunc}   For $x\in\Gamma_0$, let $L_x$ be the simple  $k\gl$-module associated to $x$.  Let $M\in\mathrm{f.d.}\,k\gl$ be indecomposable.  
\begin{itemize}
\item[(a)]  Suppose $x$ is a sink in  $\gl$.  If $M\cong L_x$ then $F^{+}_{x}M=0$.   If $M\not\cong L_x$, then $F^{+}_{x}M$ is indecomposable and $F^{-}_{x}F^{+}_{x}M\cong M$.
\item[(b)]  Suppose $x$ is a source in  $\gl$.   If $M\cong L_x$ then $F^{-}_{x}M=0$.   If $M\not\cong L_x$, then $F^{-}_{x}M$ is indecomposable and $F^{+}_{x}F^{-}_{x}M\cong M$.
\end{itemize}
\end{thm}

 Let us denote 
by $\mathfrak S$ the set of (+)-admissible sequences on $\gl$. If  $S=x_{1}¥,
\dots,x_{s}¥$ is in $\mathfrak S$,  we set
$F(S)=F^{+}_{x_{s}}\dots F^{+}_{x_{1}}:\mathrm{f.d.}\,k\gl\to
\mathrm{f.d.}\,k(\Gamma,\Lambda^{S})$.
If $S$ consists of 
distinct vertices and contains each 
vertex of the quiver, then $F(S)=\Phi^{+}$ does not depend on the 
choice of $S$ and is called 
the {\it positive Coxeter functor} ~\cite[Definition 1.2]{bgp}. For $S\in\mathfrak S$ 
we say that $S$ {\it annihilates} a module 
$M\in\mathrm{f.d.}\,k\gl$ if $F(S)M=0$.  

Replacing sinks with 
sources, one gets similar definitions of a
($-$)-admissible sequence and the negative Coxeter functor $\Phi^{-}¥$. 

We now quote 
~\cite[Definition 1.3]{bgp}.
\begin{defn}\label{intr1} A module $M\in\mathrm{f.d.}\,k\gl$ is 
preprojective if $(\Phi^{+}¥)^{m}¥M=0$ for some integer 
$m\ge0$.
\end{defn}
According to ~\cite[\S 1, Note 2]{bgp}, Definition \ref{intr1} is 
equivalent to the following.
\begin{defn}\label{intr2} A module $M\in \mathrm{f.d.}\,k\gl$ is 
 preprojective if there exists an $S\in\mathfrak S$ that 
annihilates it.
\end{defn}

\section{Lattice of (+)-admissible Sequences}\label{Section:lattice}

We begin by recalling some of the results of ~\cite{kt} needed in the sequel. 

\begin{defn}\label{mult} For $S=x_1, \dots, x_s$ in $\mathfrak{S}$ with $s\ge0$, we define the {\it length} of $S$ as $\ell(S)=s$; the {\it support} of $S$  as $\mathrm{Supp}\,S=\{v\in\Gamma_0\,|\, \exists j, 0< j\le s, \mathrm{with}\ v=x_j\}$; and for all $v\in\Gamma_0$, the {\it multiplicity} of $v$ in $S$, $m_S(v)$, as the (nonnegative) number of subscripts $j$ satisfying $x_j=v$.  A sequence $K\in\mathfrak{S}$ is  {\it complete} if $m_K(v)=1$ for all $v\in\Gamma_0$.  If $K\in\mathfrak{S}$ is  complete, $\La^K=\La$ so that if $m>0$ and $K^m$ denotes the concatenation of $m$ copies of $K$, we have $K^m\in\mathfrak{S}$. 
\end{defn} 

The following is ~\cite[Definition 1.2]{kt}.

\begin{defn}\label{def:sim}
If a sequence ${S = x_1, \dots x_{i},x_{i+1}, \dots, x_s,0<i<s}$, in
$\mathfrak S$ has the property that no edge of ${\Gamma}$
connects ${x_{i}}$ with ${x_{i+1}}$, then ${T=x_1, \dots, x_{i+1},
x_{i}, \dots x_s}$ is in ${\mathfrak S}$, and we set ${SrT}$. We
denote by ${\sim}$ the equivalence relation that is a reflexive and
transitive closure of the symmetric binary relation ${r}$.
\end{defn}

The equivalence relation $\sim$ was motivated by the fact that \ds{S \sim
T} implies ${F(S) = F(T)}$ ~\cite[Lemma 1.2, proof of part 3)]{bgp}.    

The following statement, which is ~\cite[Proposition 1.9]{kt}, gives
a canonical form in ${\mathfrak{S}}$.  

\begin{prop}\label{Prop:canonical}
Let ${S \in \mathfrak{S}}$ be nonempty.

\begin{itemize}

\item[(a)]We have ${S \sim S_1 S_2 \dots S_r}$, the concatenation of $S_1,\dots,S_r$, where, for all ${i}$,
${S_{i}}$ consists of distinct vertices, and
${\mathrm{Supp}\,{S_{i}}=\mathrm{Supp}\,{S_{i}S_{i+1}\dots S_{r}}}$. Further, if
${\mathrm{Supp}\,{S_{i}}\neq\Gamma_{0}}$ then
\newline${\mathrm{Supp}\,{S_{i+1}}\subsetneq\mathrm{Supp}\,{S_{i}}}$.

\item[(b)] Let ${T\sim T_{1}T_{2}\dots T_{q}\,}$ be a nonempty sequence
in ${\mathfrak S}$ where, for all ${j}$, ${T_{j}}$ consists of
distinct vertices, and ${\mathrm{Supp}\,{T_{j}}=\mathrm{Supp}\,{T_{j}T_{j+1}\dots
T_{q}}}$.  Then ${S\sim T}$ if and only if ${r=q}$ and
${S_{i}\sim T_{i}}$ on ${(\Gamma,\Lambda^{S_{1}\dots
S_{i-1}}),i=1,\dots,r}$.

\end{itemize}
\end{prop}

For  ${S \in
\mathfrak{S}}$, the sequence ${S_1 S_2 \dots S_r}$ of Proposition \ref{Prop:canonical}(a) is called
the \textit{canonical form} of $S$,  the integer $r$ is the {\it size} of $S$, and $S_i$ is the $i$th {\it segment} of $S$. 

\begin{rmk}\label{multcanon} In the setting of Proposition \ref{Prop:canonical}(a), if $v\in \Gamma_0$  then $v\in\mathrm{Supp}\,S_i$  if and  only if $m_S(v)\ge i$.
\end{rmk}

The sequences in $\mathfrak S$ are classified up to equivalence in terms of filters of $(\Gamma_0,\La)$.  Recall that a subset $F$ of a poset $(P,\le)$  is a {\it filter} if for all $x\in F$ and $y\in P$, $x\le y$ implies $y\in F$; a filter $F$ is {\it principal} if $F=\langle x\rangle=\{y\in P\,|\,x\le y\}$.  For a filter $F$ of $(\Gamma_0,\La)$, the {\it hull} of $F$, $H_{\La}(F)$, is the smallest filter of $(\Gamma_0,\La)$ containing $F$ and each vertex  of $\Gamma_0$ connected by an edge to a vertex in $F$.  

\begin{rmk}\label{ConnectHull} If $F$ is a filter of  $(\Gamma_0,\La)$ such that the full subgraph of $\Gamma$ determined   by $\mathrm{Supp}\,F$ is connected (for example,  if $F$ is a principal filter), then the full subgraph of $\Gamma$ determined   by $\mathrm{Supp}\,H_{\La}(F)$ is connected.
\end{rmk} 

If $X\subset \Gamma_0$, there exists an $S\in\mathfrak S$ satisfying $\mathrm{Supp}\,S=X$ if and only if $X$ is a filter of $(\Gamma_0,\La)$, and  if $\mathrm{Supp}\,S=X$ and the vertices of $S$ are distinct, then $S$ is unique up to equivalence ~\cite[Proposition 1.3]{kt}.  We now recall the classification of sequences in $\mathfrak S$ given in  ~\cite[Proposition 1.11]{kt}.

\begin{prop}\label{tails} {\rm (a)} If $S=S_{1}S_{2}\dots S_{r}\in\mathfrak S$ is
a nonempty sequence in the canonical form  then, for 
all $i$, $\mathrm{Supp}\,{S_{i}}$ 
is a filter of $(\Gamma_{0}¥,\Lambda)$ and, for $0<i<r$, 
$H_{\Lambda}(\mathrm{Supp}\,{S_{i+1}})\subset\mathrm{Supp}\,{S_{i}}$.

{\rm (b)} If $F_{1}\supset\dots\supset F_{r-1}\supset F_{r}$ is a sequence 
of nonempty
filters of $\,(\Gamma_{0}¥,\Lambda)$ satisfying
$H_{\Lambda}(F_{i+1})\subset F_{i}$ for $\,0<i<r$, then there exists a 
unique up to equivalence sequence $S_{1}S_{2}\dots S_{r}\in\mathfrak S$ in the 
canonical 
form satisfying $\mathrm{Supp}\,{S_{i}}=F_{i}$ for all $i$.
\end{prop}

The following  is ~\cite[Definition
2.1]{kt}.

\begin{defn}
If ${S,T \in \mathfrak{S}}$, we say that ${S}$ is a
{\it subsequence} of ${T}$ and write ${S\preccurlyeq T}$ if 
${T\sim SU}$ for some (+)-admissible sequence ${U}$.
\end{defn}

It was shown in \cite[Section 2]{kt} that ${\preccurlyeq}$ is a
preorder on ${\mathfrak{S}}$, and that ${S\preccurlyeq T}$ and
${T\preccurlyeq S}$ if and only if ${S \sim T}$.  Therefore the preorder $\preccurlyeq$ induces a partial order on the set of equivalences classes of sequences in $\mathfrak S$.  We often identify equivalent  (+)-admissible sequences and then say that $\preccurlyeq$ is a partial order on $\mathfrak S$.  
The next statement is a characterization of
 the preorder in terms of the canonical form.
 
\begin{prop}\label{Prop:TFAE}
If \ds{S,T\in\mathfrak S} are nonempty and if \ds{S_1 \dots S_r}, \ds{T_1 \dots
T_q} are their canonical forms, respectively, then the following are
equivalent.

\begin{itemize}
\item[(a)]\ds{S \preccurlyeq T}.
\item[(b)]\ds{r \leq q} and \ds{\mathrm{Supp}\,{S_i} \subset
\mathrm{Supp}\,{T_i}} for \ds{0 < i \leq r}.
\item[(c)] For all $v\in\Gamma_0$, $m_S(v)\le m_T(v)$.
\end{itemize}

\end{prop}

\begin{proof}
The fact that \ds{S \preccurlyeq T} if and only if \ds{r \leq q} and
\ds{S_i \preccurlyeq T_i} for \ds{0 < i \leq r} is ~\cite[Proposition
2.1(c)]{kt}.  Since ${S_i},{T_i} \in \mathfrak S$  consist
of distinct vertices,  \ds{S_i \preccurlyeq T_i} is equivalent to
\ds{\mathrm{Supp}\,{S_i} \subset \mathrm{Supp}\,{T_i}} according to ~\cite[Proposition 1.6, parts (a) and
(b)]{kt}.  Therefore (a) is
equivalent to (b).

The fact that (b) and (c) are equivalent is an immediate consequence of Remark \ref{multcanon}.
\end{proof}

\begin{cor}\label{multequiv} If \ds{S,T\in\mathfrak S}, then \ds{S \sim T} if and only if for all $v\in\Gamma_0$, $m_S(v)= m_T(v)$.
\end{cor}

The relations \ds{\sim} and  \ds{\preccurlyeq}
satisfy the left cancellation property.

\begin{prop}\label{Prop:cancellation} Let  \ds{S \in
\mathfrak{S}} and  let ${U, V}$ be (+)-admissible sequences on
\quiver{\Gamma}{\Lambda^S}.
\begin{itemize}

\item[(a)]\ds{SU \preccurlyeq SV}
if and only if \ds{U \preccurlyeq V}.
\item[(b)]\ds{SU \sim SV} if and only if
\ds{U \sim V}.

\end{itemize}
\end{prop}

\begin{proof}  Part (a) is an immediate consequence of the equivalence of parts (a) and (c) of
Proposition \ref{Prop:TFAE}, and (b) follows directly from Corollary \ref{multequiv}.
\end{proof}

To show that the poset  of equivalence classes of (+)-admissible sequences is  a lattice, we define the greatest lower and the least upper bounds, $\wedge$ and $\vee$.  

\begin{defn}\label{def:meetjoin}
Let \ds{S, T \in \mathfrak{S}} be nonempty and let
 \ds{S_{1}S_{2} \dots S_{r}}, \ds{ T_{1}T_{2} \dots T_{q}} be their canonical forms, respectively, where without loss of generality we assume that \ds{r \leq q}.  We set:
\begin{itemize}
\item[(a)] ${S \wedge T}$ to be a (+)-admissible sequence with the
canonical form \ds{R_{1}R_{2} \dots R_{r}} where \ds{\mathrm{Supp}\,{R_{i}} = \mathrm{Supp}\,{S_{i}} \cap \mathrm{Supp}\,{T_{i}}} for \ds{0< i \leq r}. 

\item[(b)] ${S \vee T}$ to be a (+)-admissible sequence with the
canonical form \ds{R_{1}R_{2} \dots R_{q}} where \ds{\mathrm{Supp}\,{R_{i}} = \mathrm{Supp}\,{S_{i}} \cup \mathrm{Supp}\,{T_{i}}} for \ds{0< i \leq r}, and  \ds{\mathrm{Supp}\,{R_{i}} = \mathrm{Supp}\,{T_{i}}} for
\ds{r < i \leq q}.
\end{itemize}

If $\emptyset$ is the empty sequence in $\mathfrak{S}$, then for all $S\in\mathfrak{S}$, we set  ${S \wedge \emptyset=\emptyset}$  and ${S \vee \emptyset=S}$.
\end{defn}

That $S \wedge T$ and $S
\vee T$ are in fact (+)-admissible sequences is contained in the proof of the
following proposition.

\begin{prop}\label{lattice}
The poset of equivalence classes of $\sim$ in \ds{\mathfrak{S}} with the partial order
\ds{\preccurlyeq} is a lattice where the operations of the greatest lower bound and the least upper bound  are \ds{\wedge} and \ds{\vee}, respectively.
\end{prop}

\begin{proof} The intersection or union of two filters is always a filter.  If $F_1,F_2$ are filters of $(\Gamma_0,\La)$, then it is straight forward that $H_{\La}(F_1\cap F_2)\subset H_{\La}(F_1)\cap H_{\La}(F_2)$ and $H_{\La}(F_1\cup F_2)= H_{\La}(F_1)\cup H_{\La}(F_2)$.  Therefore, in view of Proposition \ref{tails}, we conclude that if $S,T\in\mathfrak S$, then 
${S \wedge T}$ and ${S \vee T}$ are in \ds{\mathfrak{S}}.   It is easy to check that ${S \wedge T}$ and ${S \vee T}$  are the
greatest lower bound and the least upper bound, respectively, of $S$ and $T$.  We leave it to the reader.
\end{proof}

The following statement is a generalization of ~\cite[Lemma 1.4]{kt}

\begin{thm}\label{Theorem:meetjoin}
Let \ds{S, T \in \mathfrak{S}} be nonempty.
\begin{itemize}
 \item[(a)]  \ds{S \sim \left(S
\wedge T \right)U}, \ds{T \sim \left(S \wedge T \right)V} where $U,V$ are (+)-admissible sequences on \quiver{\Gamma}{\Lambda^{S \wedge T}} that are unique up to equivalence.
\item[(b)]\ds{\mathrm{Supp}\,U \cap \mathrm{Supp}\,V = \emptyset}.
\item[(c)]\ds{UV}, \ds{VU} are (+)-admissible sequences on \quiver{\Gamma}{\Lambda^{S
\wedge T}} and \ds{UV \sim VU}.
\item[(d)]\ds{S \vee T \sim \left(S \wedge T \right)UV \sim SV \sim TU}.

\end{itemize}
\end{thm}

\begin{proof}  (a)  This is a direct consequence of Propositions \ref{lattice} and \ref{Prop:cancellation}(b).

(b)  By (a), we have $(S \wedge T) (U \wedge V)\preccurlyeq S,T $, so Proposition \ref{lattice}  implies $(S \wedge T) (U \wedge V)\preccurlyeq S\wedge T$ whence $U \wedge V=\emptyset$.  By Definition \ref{def:meetjoin}(a) and Proposition \ref{Prop:canonical}(a), \ds{\mathrm{Supp}\,U \cap \mathrm{Supp}\,V = \emptyset}.

(c)  Since \ds{\mathrm{Supp}\,U} is a
filter, there is no arrow $u_i\to v_j$ in $\gl$ with \ds{u_i \in \mathrm{Supp}\,U}, 
\ds{v_j \in \ds{\mathrm{Supp}\,V}}, and a similar conclusion holds for \ds{\mathrm{Supp}\,V}.  Now the statement follows
immediately from (b).

(d)  By (a) and Proposition \ref{lattice}, we have 
${S \vee T \sim \left(S \wedge T \right
)U V' \sim \left(S \wedge T \right )V U'},$ 
for some $U',V'$, as well as $S,T\preccurlyeq(S \wedge T)UV\sim(S \wedge T )VU $, using (c). By Proposition \ref{lattice}, 
$${S \vee T \sim \left(S \wedge T \right
)U V' \sim \left(S \wedge T \right )V U'}\preccurlyeq(S \wedge T)UV\sim(S \wedge T )VU.$$
Applying the cancellation laws of Proposition \ref{Prop:cancellation} to the displayed formulas, we get \ds{UV' \sim VU'} and $V'\preccurlyeq V, U'\preccurlyeq U$.  By (b),  $m_{U'}(u)=m_{U}(u)$ for $u\in\mathrm{Supp}\,U$ and $m_{V'}(v)=m_{V}(v)$ for $v\in\mathrm{Supp}\,V$, so Corollary \ref{multequiv} implies  \ds{U' \sim U} and 
\ds{V' \sim V}.  Thus (d) holds. 
\end{proof}

\section{Principal admissible sequences} \label{S:principal}

We quote ~\cite[Definition 2.2]{kt}.

\begin{defn}\label{tight} A nonempty sequence $S\in\mathfrak S$ is {\it principal}
 if  its 
canonical form $S_{1}S_{2}\dots S_{r}$ satisfies 
$\mathrm{Supp}\,{S_{i}}=H_{\La}(\mathrm{Supp}\,{S_{i+1}})$ for $0<i<r$ where $\mathrm{Supp}\,{S_{r}}$ is a
principal filter. We denote by  $\mathfrak P$ the set 
of principal sequences in $\mathfrak S$. By Proposition \ref{tails}, a principal sequence is 
determined uniquely up to equivalence by its size $r$ and the set $\mathrm{Supp}\,{S_{r}}$, so we let $S_{r,x}¥$ 
denote the principal 
sequence of size $r$ with $\mathrm{Supp}\,{S_{r}}=\langle x\rangle,\ 
x\in\Gamma_{0}$. Thus $\mathfrak P=\{S_{r,x}\,|\,r\in\mathbb 
Z^{+},x\in\Gamma_{0}\}$ where $\mathbb Z^{+}$ is the set of positive integers.
\end{defn}

\begin{rmk}\label{ConnectPrinc}  It follows from Remark \ref{ConnectHull} that if $S\in\mathfrak P$, the full subgraph of $\Gamma$ determined by $\mathrm{Supp}\,S$ is connected.
\end{rmk}
We quote ~\cite[Corollary 2.3]{kt}.

\begin{prop}\label{PrincLess}  Let $S,T\in\mathfrak S$ be nonempty, let $S_1\dots S_r$ be the canonical form of $S$, and let $T=y_1,\dots,y_t$ be in $\mathfrak P$.  If $T\sim S_{q,y}$ then:
\begin{itemize}
\item[(a)]  $T\preccurlyeq S$ if and only if $q\le r$ and $y\in\mathrm{Supp}\,{S_{q}}$.
\item[(b)]  $y_t=y$.
\end{itemize}
\end{prop}

A nonempty sequence in  $\mathfrak S$ is the join of some sequences in $\mathfrak P$.

\begin{prop}\label{SeqJoinPrinc}  Let $\emptyset\ne S\in\mathfrak S$ and let $S_1\dots S_r$ be the canonical form of $S$.   Set $S_{r+1}=\emptyset$ and  $\mathrm{Supp}\, S_{r+1}=\emptyset$.
\begin{itemize}
\item[(a)]   $S=\underset{(h,v)}\bigvee S_{h,v}$ where $0<h\le r$ and, for a given $h$, $v$ runs through the set of minimal elements of $\mathrm{Supp}\, S_{h}\setminus H_{\La}(\mathrm{Supp}\, S_{h+1})$ in the partial order of $(\Gamma_0,\La)$.\vskip.05in
\item[(b)]  If  $S=T_1\vee\dots\vee T_l$ where $T_i\in \mathfrak P$ for all $i$, then for each pair $(h,v)$ described in (a), there exists an $i$ satisfying $S_{h,v}\sim T_i$.\vskip.03in
\item[(c)]  There exist $T_1,\dots,T_l\in \mathfrak P$ satisfying $S=T_1\vee\dots\vee T_l$.  If $l$ is the smallest possible and $S=U_1\vee\dots\vee U_l$ where $U_1,\dots,U_l\in \mathfrak P$, there exists a reindexing so that $T_i\sim U_i$ for all $i$.
\end{itemize}
\end{prop}
\begin{proof} (a) Proceed by induction on $r$.  If $r=1$, then $h=1$ in all pairs $(h,v)$ and $\mathrm{Supp}\, S=\mathrm{Supp}\, S_{1}$ is a filter of $(\Gamma_0,\La)$.  Since a nonempty filter is the union of the principal filters generated by its minimal elements, the statement follows from Definition \ref{def:meetjoin}(b).  Suppose now that $r>1$ and the statement holds for all nonempty sequences in $\mathfrak S$ of size $<r$. By the induction hypothesis, $S_2\dots S_r=\underset{(h,v)}\bigvee S_{h-1,v}$ where $1<h\le r$.  It follows from Definition \ref{def:meetjoin}(b) that $S=(\underset{(1,v)}\bigvee S_{1,v})\bigvee(\underset{(h,v)}\bigvee S_{h,v})$ where $1<h\le r$.  If $u\in\mathrm{Supp}\, S_{1}\setminus H_{\La}(\mathrm{Supp}\, S_{2})$ satisfies $u<v$ in the poset $(\Gamma_0,\La)$, then $S_{1,v}\preccurlyeq S_{1,u}$ in $\mathfrak S$.  Therefore $\underset{(1,v)}\bigvee S_{1,v}=\underset{(1,u)}\bigvee S_{1,u}$ where $u$ runs through the set of minimal elements of $\mathrm{Supp}\, S_{1}\setminus H_{\La}(\mathrm{Supp}\, S_{2})$.  The proof of (a) is complete.

(b) Suppose $S=T_1\vee\dots\vee T_l$ where $T_i\in \mathfrak P$ for all $i$.  Since  $v\in\mathrm{Supp}\, S_{h}$ for each $(h,v)$, Definition \ref{def:meetjoin}(b) implies that there exists an $i$ such that the canonical form $W_1\dots W_q$ of $T_i$ satisfies $v\in\mathrm{Supp}\, W_{h}$.  Hence $h\le q$ and Proposition \ref{PrincLess}(a) says that $S_{h,v}\preccurlyeq T_i$.  Since $T_i\in \mathfrak P$, we have $T_i\sim S_{p,u}$ for some $p>0$ and $u\in\Gamma_0$, whence $u\in\mathrm{Supp}\, S_{p}$.  Therefore there exists a pair $(j,w)$, where $w$ is a minimal element of $\mathrm{Supp}\, S_{j}\setminus H_{\La}(\mathrm{Supp}\, S_{j+1})$, such that the canonical form $X_1\dots X_j$ of $S_{j,w}$ satisfies $u\in\mathrm{Supp}\, X_{p}$.  Then $p\le j$ and $T_i\preccurlyeq S_{j,w}$ whence $S_{h,v}\preccurlyeq S_{j,w}$.  By Proposition \ref{PrincLess}(a), $h\le j$ and $v\in\mathrm{Supp}\, X_{h}$.  If $h<j$, then $v\in\mathrm{Supp}\, X_{h}=H_{\La}(\mathrm{Supp}\, X_{h+1})\subset H_{\La}(\mathrm{Supp}\, S_{h+1})$, which contradicts the conditions imposed on the pair $(h,v)$ in (a).  Therefore we must have $h=j$.  Since  $S_{h,v}\preccurlyeq S_{h,w}$, then $\langle v\rangle\subset \langle w\rangle$ and $w\le v$.  The latter implies $w=v$ because $v,w$ are minimal elements of $\mathrm{Supp}\, S_{h}\setminus H_{\La}(\mathrm{Supp}\, S_{h+1})$.  It follows that $T_i\sim S_{h,v}$.

(c)  The statement is a consequence of (a) and (b).
\end{proof}

We quote  ~\cite[Definition 3.1]{kt}.

\begin{defn}If $S\in \mathfrak S$ annihilates a $k\gl$-module $M$, but no 
proper subsequence of $S$ annihilates $M$, we call $S$ a {\it                            
shortest} sequence annihilating $M$.
\end{defn}

The following statement is ~\cite[Theorems 3.1 and 3.5, Corollary 3.6(b)]{kt}.

\begin{thm}\label{shrtstsq} Let $M$ be an
indecomposable preprojective $k\gl$-module.
\begin{itemize}
\item[(a)]  There exists
a unique up to equivalence shortest  sequence
$S_M\in\mathfrak S$ annihilating $M$.  
\item[(b)]  If  $N$ is an
indecomposable preprojective $k\gl$-module, then $S_N\sim S_M$ if and only if $N\cong M$. 
\item[(c)]  $S_M\in\mathfrak P$. 
\item[(d)]  If $S_M=x_1\dots x_s$, then $M\cong F^-_{x_1}\dots F^-_{x_{s-1}}(L_{x_s})$ where $L_{x_s}$ is the simple projective \newline$k(\Gamma,\sigma_{x_{s-1}}\dots\sigma_{x_{1}}\La)$-module associated with $x_s\in\Gamma_0$.
\end{itemize}
\end{thm}

We  now drop the assumption of indecomposability of $M$ in Theorem \ref{shrtstsq}.

\begin{thm}\label{shrtst} Let $M$ be a preprojective $k\gl$-module.
\begin{itemize}
\item[(a)]  There exists
a unique up to equivalence shortest  sequence
$S_M\in\mathfrak S$ annihilating $M$.  If $M\cong M_1^{m_1}\oplus\dots\oplus M_t^{m_t}$ where the $M_i$'s are nonisomorphic indecomposable $k\gl$-modules and $m_i>0$ for all $i$, then each $M_i$ is preprojective and $S_M=S_{M_1}\vee\dots\vee S_{M_t}$. 
\item[(b)]  If $L$ is a preprojective $k\gl$-module, then $S_L\preccurlyeq S_M$ if and only if for each indecomposable direct summand $X$ of $L$, there exists an indecomposable direct summand $Y$ of $M$ satisfying  $S_X\preccurlyeq S_{Y}$.
\end{itemize}
\end{thm}
\begin{proof} (a) If $M=0$ then $S_M=\emptyset$.  If $M\ne0$, then $M\cong M_1^{m_1}\oplus\dots\oplus M_t^{m_t}$ as indicated in the statement.  Since every reflection functor is additive, each $M_i$ is  preprojective.  By Theorem \ref{shrtstsq}(a), a sequence $S\in\mathfrak S$ annihilates $M$ if and only if $S_{M_i}\preccurlyeq S$ for all $i$.  Since $\mathfrak S$ is a lattice by Proposition \ref{lattice}, we have $S_M=S_{M_1}\vee\dots\vee S_{M_t}$.  Alternatively, we note that the proof of Theorem \ref{shrtstsq}(a), ~\cite[pp. 394-395]{kt}, does not actually use the indecomposability of $M$ and works for any nonzero preprojective $M$.

(b) The statement is trivial if either $L$ or $M$ is zero.  Assuming $L, M$ are nonzero, we get $ L\cong L_1^{l_1}\oplus\dots\oplus L_s^{l_s}$ and $ M\cong M_1^{m_1}\oplus\dots\oplus M_t^{m_t}$ as in (a).  For the sufficiency, suppose that for each $i$ there exists a $j$ satisfying $S_{L_i}\preccurlyeq S_{M_j}$.  By (a), $S_{L_i}\preccurlyeq S_M$ whence $S_L=S_{L_1}\vee\dots\vee S_{L_s}\preccurlyeq S_M$ because $\mathfrak S$ is a lattice according to Proposition \ref{lattice}.
For the necessity, let $T_1\dots T_q$ be the canonical form of $S_M$ and let $X=L_i$. By Theorem \ref{shrtstsq}(c),  $S_X\in\mathfrak P$ whence $S_X\sim S_{r,x}$ where $r>0$ and $x\in\Gamma_0$. It is clear that $S_X\preccurlyeq S_L$, so $S_L \preccurlyeq S_M$ implies $S_X\preccurlyeq S_M$.      By Proposition \ref{PrincLess}(a),  $r\le q$ and $x\in\mathrm{Supp}\, T_r$.  By Definition \ref{def:meetjoin}(b), $\mathrm{Supp}\, T_r$ is the union of $r$th segments of some of the sequences $S_{M_1},\dots,S_{M_t}$.  Hence, for some $j$, the canonical form of $S_{M_j}$ is $U_1\dots U_p$ where $r\le p$ and $x\in\mathrm{Supp}\, U_r$.  By Proposition \ref{PrincLess}(a), $S_X\preccurlyeq S_{M_j}$.
\end{proof}

\begin{rmk} Part (b) of Theorem \ref{shrtstsq} is false without the assumption that both $M$  and $N$ are indecomposable.  For example, if $M$ is indecomposable and $N=L\oplus M$ where $L$ is indecomposable preprojective with $S_L\preccurlyeq S_M$, then $S_M=S_N$ but $M\not\cong N$.
\end{rmk}

Since $\mathfrak P$ is a subset of $\mathfrak S$, the partial order $\preccurlyeq$ on the set of equivalence classes of $\sim$ in $\mathfrak S$ induces a partial order on the set of equivalence classes of $\sim$ in $\mathfrak P$.  Identifying equivalent sequences in $\mathfrak P$, we often say that  $\preccurlyeq$ is a partial order on $\mathfrak P$.  The poset structure of $\mathfrak P$ carries a lot of  information about the preprojective component of the Auslander-Reiten quiver of $k(\Gamma,\La)$.  We now recall some definitions and facts from ~\cite {ars,r}. 

Let ${\mathbb N}$ be the set of nonnegative integers.  The 
translation quiver ${\mathbb N}(\Gamma,\Lambda^{op})$ of the opposite quiver of 
$\gl$ is defined as follows.  The set of vertices of ${\mathbb 
N}(\Gamma,\Lambda^{op}¥)$ is ${\mathbb N}\times\Gamma_{0}¥$, and each 
arrow $a:u\to v$ of $\gl$ gives rise to two series of arrows, 
$(n,a^{\circ}¥):(n,v)\to(n,u)$ and 
$(n,a^{\circ}¥)':(n,u)\to(n+1,v)$.  By 
construction, ${\mathbb 
N}(\Gamma,\Lambda^{op}¥)$ is a locally finite quiver without oriented 
cycles, so ${\mathbb N}\times\Gamma_{0}$ is a poset as explained earlier.
 
 Let $X\in\mathrm{f.d.}\,k\gl$ be 
 indecomposable and let 
 $[X]$ be the isomorphism class of $X$. If $Y\in\text{f.d.}\,k\gl$ is 
 indecomposable, a path of length $m>0$ from $X$ to $Y$ is a 
 sequence of 
 nonzero nonisomorphisms $X=A_{0}¥\to\dots\to A_{m}¥=Y,$ where
 $A_{i}¥\in\text{f.d.}\,k\gl$ is
 indecomposable for all $i$. By definition, there 
  exists a path of 
 length zero from $X$ to $X$.  One writes $[X]\prec [Y]$ if there
  exists a path of positive length from $X$ to $Y$. 
  
  The preprojective component of $\gl$,
 $\tilde{\mathscr P}\gl$, is a locally finite connected 
 quiver whose set of vertices, $\tilde{\mathscr P}\gl_{0}¥$, consists 
 of the isomorphism classes of 
 indecomposable 
 preprojective $k\gl$-modules, and the number of arrows 
 $[X]\to[Y]$ is the $k$-dimension of the vector space Irr$\,(X,Y)$ of irreducible 
 maps $X\to Y$.  If 
 $X,Y$ are indecomposable where $Y$ is preprojective, and if $X=A_{0}¥\to\dots\to 
 A_{m}¥=Y,\ m>0,$ is a path from $X$ to $Y$, then $[X]\ne[Y]$ and
 $A_{i}¥$ is preprojective for all $i$.  It follows that the reflexive closure 
 $\preccurlyeq$ 
 of the transitive binary relation $\prec$ is a partial order on 
 $\tilde{\mathscr P}\gl_{0}¥$.  
 Moreover, $[X]\prec [Y]$ if and only if there is a finite 
 sequence of 
 irreducible morphisms $X=B_{0}¥\to\dots\to B_{n}¥=Y$, where $n>0$ 
 and $B_{j}¥$ is
 indecomposable preprojective for all $j$.  
 
 \begin{thm}\label{translquiv}
 \begin{itemize}
\item[(a)] The map $\psi:\mathfrak P\to{\mathbb N}\times\Gamma_{0}$ given by 
$\psi(S_{r,x}¥)=(r-1,x)$ 
is an isomorphism of posets.
\item[(b)] Consider the map 
 $\phi:\tilde{\mathscr P}\gl\to{\mathbb N}(\Gamma,\Lambda^{op})$ 
 defined on the vertices by 
$\phi([L])=(\nu,x)=(\nu(L),x(L))$, where $x$ is the vertex of $\gl$ 
 associated with the indecomposable projective module 
 $(\Phi^{+}¥)^{\nu}¥L\cong(\DTr)^\nu L$, and defined on the arrows in a natural way 
 ~\cite[VIII Proposition 1.15]{ars}. Given an $[M]\in\tilde{\mathscr P}\gl_{0}$, the map $\phi$ 
 induces a bijection between the set of paths in $\tilde{\mathscr P}\gl$ 
 ending at $[M]$ and the set of paths in ${\mathbb N}(\Gamma,\Lambda^{op})$ 
 ending at $\phi([M])$. 
  \item[(c)] The map $\chi:\tilde{\mathscr P}\gl_{0}¥\to\mathfrak P$ given by 
$[L]\mapsto S_{L}¥$ is an injective morphism of posets.
\item[(d)] If $\Gamma$ is not a Dynkin diagram of the type A, D, or E, the maps $\phi$ and $\chi$  are isomorphisms.
\end{itemize} 
 \end{thm}
 \begin{proof} (a) This is ~\cite[Theorem 2.5(a)]{kt}.
 
 (b) This is ~\cite[Proposition 3.7(d)]{kt}.
 
 (c) This is ~\cite[Corollary 3.8(a)]{kt}.
 
 (d) This is ~\cite[Proposition 3.7(b) and Corollary 3.8(c)]{kt}.
 \end{proof}

 We finish this section with two results that play a crucial role in Section \ref{weylgroup}.
 
 \begin{prop}\label{princseq} If $S_{r,x}\sim S=x_1,\dots,x_s,\ s>1$, then $T=x_2,\dots,x_s$ is a principal (+)-admissible sequence on $(\Gamma,\sigma_{x_1}\La)$.  If $S_1\dots S_r$ and $T_1\dots T_q$ are the canonical forms of $S$ and $T$, respectively, then $\mathrm{Supp}\, T_{q}$ is the principal filter of $(\Gamma_0,\sigma_{x_1}\La)$ generated by $x$, and we have:
 \begin{itemize}
 \item[(a)] If $x_1=x$, then $q=r-1$ and $\mathrm{Supp}\, T_i=\mathrm{Supp}\, S_i$ for $0<i<r$. 
 \item[(b)]  If $x_1\not=x$, then $q=r$, $\mathrm{Supp}\, T_i=\mathrm{Supp}\, S_i$ for $i\ne m_S(x_1)$, and $\mathrm{Supp}\, T_i=\mathrm{Supp}\, S_i\setminus\{x_1\}$ for $i= m_S(x_1)$. 
 \end{itemize}
 \end{prop} 
 \begin{proof}  Without loss of generality, we may assume that $\Gamma$ is not a Dynkin diagram of the type A, D, or E.  For if it is, there must be at least one arrow in $\gl$ because $\Gamma$ is a connected graph with more than one vertex.  We double the arrow preserving its direction.  The new graph is  no longer a Dynkin diagram, but the new quiver has the same sets $\mathfrak{P}$ and ${\mathfrak S}$ as  the original quiver had.  
 
 By Theorem \ref{translquiv}(d), the map $\chi$ of Theorem \ref{translquiv}(c)  is an  isomorphism.  Hence $S\sim S_M$ for some indecomposable preprojective $k\gl$-module $M$, and $T=S_{F^+_{x_1}M}$ where, by Theorem \ref{PosReflFunc}(a), $F^+_{x_1}M$  is an indecomposable preprojective $k(\Gamma,\sigma_{x_1}\La)$-module because $s>1$.  By Theorem \ref{shrtstsq}(c), $T$ is a  principal (+)-admissible sequence on $(\Gamma,\sigma_{x_1}\La)$.  Since $S_{r,x}\sim S$, Proposition \ref{PrincLess}(b) says that $x_s=x$ and $\mathrm{Supp}\, T_{q}$ is the principal filter of $(\Gamma_0,\sigma_{x_1}\La)$ generated by $x$. We also have $m_T(v)=m_S(v)$  if $x_1\ne v\in\Gamma_0$, and $m_T(x_1)=m_S(x_1)-1$.  Comparing the  multiplicities of vertices in $S$ and $T$, using  Remark \ref{multcanon}, and taking into account that $\mathrm{Supp}\, S_r=\{x\}$ if $x_1=x$, we see that (a) and (b) hold.
 \end{proof}
 
The statement of Proposition \ref{princseq} does not involve representation theory, and we have a purely combinatorial proof that is longer and more technical than the  one given above.  As we noted in the proof, if $S\sim S_M$ where $M$ is an indecomposable preprojective $k\gl$-module, then $T\sim S_{F^+_{x_1}M}$ where  $F^+_{x_1}M$  is an indecomposable preprojective $k(\Gamma,\sigma_{x_1}\La)$-module.  Since an indecomposable preprojective module is uniquely up to isomorphism determined by the shortest (+)-admissible sequence that annihilates it (Theorem \ref{shrtstsq}(b)), the explicit computation of the canonical form of $T$ in terms of the canonical form of $S$  allows us to think of a positive reflection functor as operating on  principal (+)-admissible sequences instead of indecomposable preprojective  modules.  In particular, knowing the  pair $(r, x)$, which determines the location of $M$ in the preprojective component of $\gl$, we compute the pair $(q,x)$ that determines the location of $F^+_{x_1}M$ in the preprojective component of $(\Gamma,\sigma_{x_1}\La)$ (see Theorem \ref{translquiv}).
 
 \begin{prop}\label{fullsubposet} If $M,N$ are
indecomposable preprojective $k\gl$-modules, then $[M]\preccurlyeq[N]$  in $\tilde{\mathscr P}\gl_{0}$ if and only if $S_M\preccurlyeq S_N$ in $\mathfrak P$.
\end{prop}
\begin{proof} The necessity is an immediate consequence of Theorem \ref{translquiv}(c).  We now assume that $S_M\preccurlyeq S_N$ and show that $[M]\preccurlyeq[N]$.  If $S_N\preccurlyeq S_M$, then $S_M\sim S_N$ so Theorem \ref{shrtstsq}(b) implies $M\cong N$ and $[M]=[N]$;  in particular, $[M]\preccurlyeq[N]$.  Suppose now that $S_N\not\preccurlyeq S_M$ where $S_M\sim S_{p,u}$ and $S_N\sim S_{q,v}$ for some $p,q>0$ and $u,v\in\Gamma_0$.  By Theorem \ref{translquiv}(a), $(p-1,u)<(q-1,v)$ in ${\mathbb N}\times\Gamma_{0}$ whence there is a path $(p-1,u)\to(q-1,v)$ of positive length in ${\mathbb N}(\Gamma,\Lambda^{op})$.  By Theorem \ref{translquiv}(b), there is a path $[M]\to[N]$ of positive length in $\tilde{\mathscr P}\gl$, i.e., $[M]\preccurlyeq[N]$.
\end{proof} 

\section{Reduced words in the Weyl group }\label{weylgroup}

Let $A=(a_{ij})$ be an indecomposable symmetric generalized $n\times n$ Cartan matrix (see ~\cite{kac1990}),  i.e.,  $A$ is a symmetric integral matrix with $a_{ii}=2$ for all $i$ and $a_{ij}\le0$ for $i\ne j$ that is not conjugate under a permutation matrix to a block-diagonal matrix $\left(\begin{matrix}A_1&0\\0&A_2\end{matrix}\right)$.  For the rest of the paper, we fix the matrix $A$ and assume that in the graph $\Gamma=(\Gamma_0,\Gamma_1)$ we have $\Gamma_0=\{1,\dots,n\}$ and, for all $i\ne j$, $-a_{ij}$ edges connect vertices $i$ and $j$.  To any finite connected graph  without loops, there corresponds a unique up to conjugation by a permutation matrix indecomposable symmetric generalized  Cartan matrix.  Therefore our assumptions impose no additional restrictions on $\Gamma$.  We identify the  root lattice $Q$ associated with $A$ with  the free abelian group ${\mathbb Z}^n$ by  identifying  the simple roots $\alpha_1,\dots,\alpha_n$ of $Q$ with the standard basis vectors $e_1,\dots,e_n$ of ${\mathbb Z}^n$, and we think of the latter vectors as indexed by the vertices of $\Gamma$.  Then the simple reflections on $Q$ identify with the reflections $\sigma_1,\dots,\sigma_n$ on ${\mathbb Z}^n$ given by $\sigma_i(e_j)=e_j-a_{ij}e_i$ for all $i,j$, and the Weyl group $\mathcal W$ is the subgroup of $GL({\mathbb Z}^n)$ generated by $\sigma_1,\dots,\sigma_n$. In view of the above identification, the terms  \lq\lq root lattice" and \lq\lq Weyl group" make sense for the graph  $\Gamma$ ~\cite[Definition 2.1]{bgp}.

We quote ~\cite[Theorem 1.1]{bgp}.
\begin{thm}\label{PosRefl} If $x$ is a sink (respectively, source) in $\gl$ and $M\in\mathrm{f.d.}\,k\gl$ is indecomposable and not simple projective (respectively, injective),  then $\dim F^+_{x}M=\sigma_x(\dim M)$ (respectively, $\dim F^-_{x}M=\sigma_x(\dim M)$).
\end{thm}

\begin{defn}If $S=x_{1},\dots,x_{s}$ is in $\mathfrak S$, we
set $w(S)=\sigma_{x_{s}}\dots  \sigma_{x_{1}}$ and say that $w(S)$ is the word in the Weyl group $\mathcal W$ associated to $S$.  If no edge connects vertices $i$ and $j$, then $\sigma_i\sigma_j=\sigma_j\sigma_i$ so that $S\sim T$ implies $w(S)=w(T)$.
\end{defn}

To illustrate the utility of words associated to sequences in $\mathfrak S$, we begin with an elementary proof of the following well known fact (see ~\cite[VIII Corollary 2.3] {ars}).

\begin{prop} \label{reduced1} Let $M,N\in\mathrm{f.d.}\,k\gl$ be indecomposable.  If $M$ is preprojective and $\dim M=\dim N$, then $M\cong N$.
\end{prop}
\begin{proof}  If $S_M=x_1,\dots,x_s$ and $T=x_1,\dots,x_{s-1}$, then $F(S_M)M=0$ but $F(T)M\not=0$.  Using Theorems \ref{PosReflFunc}(a) and \ref{PosRefl},  we obtain $w(S_M)(\dim M)=w(S_M)(\dim N)<0$.  Using the  same theorems, we get $F(S_M)N=0$ whence $N$ is preprojective and $S_N\preccurlyeq S_M$. Since  $N$  is preprojective, by the symmetry we get $S_M\preccurlyeq S_N$ whence $S_M\sim S_N$.  By Theorem \ref{shrtstsq}(b), $M\cong N$.
\end{proof}

Recall (see ~\cite{bourbLie4-6}) that for $w\in\mathcal W$,  the {\it length} of $w$, $\ell(w)$, is the smallest integer $l\ge0$ such that $w$ is the product of $l$ simple reflections, and a word $w=\sigma_{y_{t}}\dots  \sigma_{y_{1}}$ in $\mathcal W$ is {\it reduced} if $\ell(w)=t$. 

\begin{rmk}\label{SubseqReduced} If $S\preccurlyeq T$ in $\mathfrak S$ where  $w(T)$ is reduced, then $T\sim SU$ for some $U$, and $w(S),w(U)$ are reduced, as follows from  $w(T)=w(U)w(S)$.
\end{rmk}

Recall that if  $v_1,\dots,v_n$ are distinct vertices of $\Gamma$, then $c=\sigma_{v_n}\dots\sigma_{v_1}$ is a {\it Coxeter element} of $\mathcal W$ (a {\it Coxeter transformation} in ~\cite[Defintion 2.3]{bgp}); $c$ depends on the choice of the permutation $v_1,\dots,v_n$ of the vertices $1,\dots,n$.

We examine the words in the Weyl group associated to preprojective modules.

\begin{thm}\label{reducedword} Let $M$ be a preprojective $k\gl$-module.
\begin{itemize}
\item[(a)]  The word $w(S_M)\in \mathcal W$ is reduced.
\item[(b)]  If $M$ is
indecomposable and $N$ is an
indecomposable preprojective $k\gl$-module, the following are equivalent.
\begin{enumerate}
\item[(i)] $M\cong N$.
\item[(ii)] $S_M\sim S_N$.
\item[(iii)] $w(S_M)=w(S_N)$.
\end{enumerate}
 
\end{itemize}
\end{thm}
\begin{proof} (a) If $M=0$ the statement is trivial.  If $M\ne0$, let $S_M=x_1,\dots,x_s$ and proceed by induction on $s>0$.  The case $s=1$ is clear, so suppose $s>1$ and the statement holds for all orientations $\Theta$ on $\Gamma$ without oriented cycles and all preprojective $k(\Gamma,\Theta)$-modules $N$ for which $S_N=y_1,\dots,y_t$ satisfies $t<s$.  Since $s>1$, $F^+_{x_1}M\ne0$ is a preprojective  $k(\Gamma,\sigma_{x_1}\La)$-module, and $S_{F^+_{x_1}M}=x_2,\dots,x_s$.  By the induction hypothesis, the word $u=\sigma_{x_s}\dots\sigma_{x_2}$ in $\mathcal W$ is reduced.  Assume, to the contrary, that the word $u\sigma_{x_1}=\sigma_{x_s}\dots\sigma_{x_2}\sigma_{x_1}$ is not reduced.  Then $\ell(u\sigma_{x_1})\le\ell(u)$ and, since $\mathcal W$ is a Coxeter group ~\cite[Proposition 3.13]{kac1990}, we must have $\ell(u\sigma_{x_1})<\ell(u)$ ~\cite[Ch. IV, Proposition 1.5.4]{bourbLie4-6}.  By ~\cite[Lemma 3.11,  part a)]{kac1990}, $\sigma_{x_s}\dots\sigma_{x_2}(e_{x_1})<0$ where $e_{x_1}$ is the simple root associated to the vertex $x_1$, whence $F(S_{F^+_{x_1}M})L_{x_1}=0$ where $L_{x_1}$ is the simple $k(\Gamma,\sigma_{x_1}\La)$-module associated to $x_1$, as follows from Theorems \ref{PosReflFunc}(a) and \ref{PosRefl}.    Since $x_1$ is a sink in $\gl$, it is a source in $(\Gamma,\sigma_{x_1}\La)$ so  $L_{x_1}$ is a simple injective and a preprojective $k(\Gamma,\sigma_{x_1}\La)$-module.  In particular, $[L_{x_1}]$ is a sink in $\tilde{\mathscr P}(\Gamma,\sigma_{x_1}\La)$ and, hence, a maximal element of the poset $\tilde{\mathscr P}(\Gamma,\sigma_{x_1}\La)_0$.  On the other hand, $S_{L_{x_1}}\preccurlyeq S_{F^+_{x_1}M}$ whence, by Theorem \ref{shrtst}(b),  $S_{L_{x_1}}\preccurlyeq S_{N}$ for some indecomposable direct summand $N$ of $F^+_{x_1}M$ and, by Proposition \ref{fullsubposet}, $[L_{x_1}]\preccurlyeq [N]$ in $\tilde{\mathscr P}(\Gamma,\sigma_{x_1}\La)_0$ .  Since  $[L_{x_1}]$ is a maximal element, we have $[L_{x_1}]=[N]$ whence $L_{x_1}\cong N$, in contradiction with the fact that the simple module associated to  a vertex that is a source is not a direct summand of a module that belongs to the image of the positive reflection functor associated to the vertex, as follows from Theorem \ref{PosReflFunc}.  Thus $w(S_M)$ is a reduced word.

(b) By Theorem \ref{shrtstsq}(b), (i) is equivalent to (ii).  It is clear that (ii)$\implies$(iii).  To prove (iii)$\implies$(ii), suppose $w(S_M)=w(S_N)$.  In view of parts (a) and (b) of Theorem \ref{Theorem:meetjoin}, we have $S_M\sim (S_M\wedge S_N)U$ and $S_N\sim (S_M\wedge S_N)V$ where $U,V$ are (+)-admissible sequences on $(\Gamma,\La^{S_M\wedge S_N})$ satisfying $\mathrm{Supp}\,U\cap\mathrm{Supp}\,V=\emptyset$.  If both $U$  and $V$ are empty, then $S_M\sim S_N$.  If not, then, say, $U=u_1,\dots,u_p$ with $p>0$.  We obtain $w(U)w(S_M\wedge S_N)=w(V)w(S_M\wedge S_N)$ whence $w(U)=w(V)$.  By (a) and Remark \ref{SubseqReduced}, the word $w(U)=\sigma_{u_p}\dots\sigma_{u_1}$ is reduced, so  ~\cite[Lemma 3.11, part b)]{kac1990} says that $w(U)(e_{u_1})<0$.  On the other hand, $w(V)(e_{u_1})$ is a root whose $u_1$-coordinate is the same as that of $e_{u_1}$, namely, is equal to $1$, because $\mathrm{Supp}\,U\cap\mathrm{Supp}\,V=\emptyset$. Hence $w(V)(e_{u_1})>0$, which contradicts  $w(U)=w(V)$.
\end{proof}

Zelevinsky suggested the following statement.

\begin{cor}\label{zel}  Let $S=x_1,\dots,x_s,\ s>0,$ be in $\mathfrak S$, and set $M(S)=F^-_{x_1}F^-_{x_2}\dots  F^-_{x_{s-1}}(L_{x_s})$, where $L_{x_s}$ is the simple projective $k(\Gamma,\sigma_{s-1}\dots\sigma_{x_1}\La)$-module associated to $x_s\in\Gamma_0$.
\begin{itemize}
\item[(a)]  If the word $w(S)\in\mathcal W$ is reduced, then $M(S)$ is an indecomposable module in $\tilde{\mathscr P}$.
\item[(b)]  If $M\in\tilde{\mathscr P}$ is indecomposable, then $M\cong M(S)$ for some sequence $S\in\mathfrak S$ where $\ell(S)>0$ and the word $w(S)$ is reduced.
\end{itemize}
\end{cor}
\begin{proof} (a) Since $w(S)$ is reduced, ~\cite[Lemma 3.10]{kac1990} implies that for $0<i<s,\ \sigma_{x_i}\dots\sigma_{x_{s-1}}(e_{x_s})>0$.  By Theorems \ref{PosReflFunc}(a) and \ref{PosRefl}, $M(S)$ is an indecomposable $k\gl$-module and $F(S)(M(S))=0$.  Hence $M(S)\in\tilde{\mathscr P}$.

(b) By Theorems \ref{shrtstsq}(d) and \ref{reducedword}(a), $M\cong M(S_M)$ and $w(S_M)$ is reduced.
\end{proof}

\begin{thm}\label{princseqreducedwrd}  If $S=x_1,\dots,x_s$, $s>0,$ is in $\mathfrak P$, the following are equivalent.
\begin{itemize}
\item[(a)]   There exists an  indecomposable preprojective $k\gl$-module $M$ satisfying $S\sim S_M$.
\item[(b)]   The word $w(S)\in \mathcal W$ is reduced
\item[(c)]  For $0<i<s,\ \sigma_{x_i}\dots\sigma_{x_{s-1}}(e_{x_s})>0$.
\end{itemize} 
\end{thm}
\begin{proof} (a)$\implies$(b)\hskip.05in  This is Theorem \ref{reducedword}(a).

(b)$\implies$(c)\hskip.05in  This is addressed in the proof of Corollary \ref{zel}(a).

(c)$\implies$(a)\hskip.05in   Set $M= M(S)=F^-_{x_1}F^-_{x_2}\dots  F^-_{x_{s-1}}(L_{x_s})$ where $L_{x_s}$ is the simple projective \newline$k(\Gamma,\sigma_{s-1}\dots\sigma_{x_1}\La)$-module associated to $x_s\in\Gamma_0$.  By Corollary \ref{zel}(a), $M$ is indecomposable preprojective.  
To show $S\sim S_M$, proceed by induction  on $s$.  The case $s=1$ is trivial, so suppose $s>1$ and the statement holds for all principal (+)-admissible sequences  of length $<s$ on all quivers $(\Gamma,\Theta)$ without oriented cycles.  

Set $N=F^-_{x_2}\dots  F^-_{x_{s-1}}(L_{x_s})$. By Proposition \ref{princseq}, $T=x_2,\dots,x_s$  is a principal (+)-admissible sequence of length $s-1$ on $(\Gamma,\sigma_{x_1}\La)$, and the same as above argument shows that $N$ is an indecomposable preprojective  $k(\Gamma,\sigma_{x_1}\La)$-module.  By the induction hypothesis, $T\sim S_N$.

It is clear that $S$ annihilates $M$, so $S_M\preccurlyeq S$ whence $S\sim S_MU$ for some $U$.   

If $x_1\in\mathrm{Supp}\, S_M$, then $S_M\sim y_1,\dots,y_t$ where $t\le s$ and $y_1=x_1$ because $x_1$ is a sink in $\gl$.  Then, using Theorem \ref{PosReflFunc}(b), we get $0=F^+_{y_t}\dots  F^+_{y_{1}}(M)=F^+_{y_t}\dots  F^+_{x_{1}}(F^-_{x_1}N)\cong F^+_{y_t}\dots  F^+_{y_{2}}(N)$ whence $ y_2,\dots,y_t$ is a (+)-admissible sequence on $(\Gamma,\sigma_{x_1}\La)$ that annihilates $N$.  Then $S_N\sim x_2,\dots,x_s\preccurlyeq y_2,\dots,y_t$ so that  $s\le t$, whence $s=t,\ U=\emptyset$, and $S\sim S_M$.

If $x_1\not\in\mathrm{Supp}\, S_M$, then $x_1\in\mathrm{Supp}\, U$ and $x_1$ is a sink in $(\Gamma,\La^{S_M})$ because $\mathrm{Supp}\, S_M$, being a filter of $(\Gamma_0,\La)$, contains  no $v\in\Gamma_0$ satisfying $v\le x_1$.  By ~\cite[Lemma 1.7]{kt}, for all $v\in\mathrm{Supp}\, S_M$, no arrow connects $v$ and $x_1$ whence $S_Mx_1\sim x_1S_M$ on $\gl$.  Therefore $S_M$ is a (+)-admissible sequence on $(\Gamma,\sigma_{x_1}\La)$ and we have $0=F^+_{x_1}(F(S_M)M)=F(S_M)(F^+_{x_1}M)=F(S_M)(F^+_{x_{1}}F^-_{x_1}N)\cong F(S_M)N$.  Hence $S_N\preccurlyeq S_M$ so that $s-1\le\ell(S_M)$, which implies $s-1=\ell(S_M)$ and $S\sim S_Mx_1\sim x_1S_M$.  Then the full subgraph of $\Gamma$ determined by $\mathrm{Supp}\, S$ is disconnected, which contradicts Remark \ref{ConnectPrinc}.
\end{proof}

\begin{thm}\label{seqreducedwrd} For all $S\in\mathfrak S$, the following are equivalent. 
\begin{itemize}
\item[(a)]  There exists a preprojective $k\gl$-module $M$ satisfying $S\sim S_M$.
\item[(b)]  The word $w(S)\in \mathcal W$ is reduced.
\end{itemize}
If $S\ne\emptyset$, let $S=T_1\vee\dots\vee T_l$ where,  for all $i$, $T_i\in\mathfrak P$ and $l$ is the smallest possible (see Proposition \ref{SeqJoinPrinc}).  Then either of (a), (b) is equivalent to the following condition.
\begin{itemize}
\item[(c)]   For $0<i\le l$, the word $w(T_i)\in \mathcal W$ is reduced.
\end{itemize}
\end{thm}
\begin{proof}  The case $S=\emptyset$ is clear, so let $S\ne\emptyset$.

(a)$\implies$(b)\hskip.05in  This is Theorem \ref{reducedword}(a).

(b)$\implies$(c)\hskip.05in  Since $T_i\preccurlyeq S$, the statement follows from Remark \ref{SubseqReduced}.

(c)$\implies$(a)\hskip.05in  By Theorem \ref{princseqreducedwrd}, $T_{i}\sim S_{M_i}$ for some indecomposable preprojective $k\gl$-module $M_i$.  Since $l$ is the smallest possible, for $i\ne j$, we have $T_i\not\sim T_j$ so that $M_i\not\cong M_j$ by Theorem \ref{shrtstsq}(b).  By Theorem \ref{shrtst}(a), $S\sim S_M$ where $M=M_1\oplus\dots\oplus M_l$.
\end{proof}

\begin{ex}  Given a graph $\Gamma$ and a reduced word $w\in\mathcal W$, it may be impossible to find an orientation $\La$ and a (+)-admissible sequence $S$ of length $\ell(w)$ on $\gl$ satisfying $w=w(S)$. 

For example, if $\displaystyle{\Gamma}=A_4$:
\[ \xymatrix{ x_1 \ar@{-}[r]^{a}  & x_2 \ar@{-}[r]  & x_3 \ar@{-}[r]^{b} &x_4},\] 
then $w=\sigma_{x_2}\sigma_{x_3}\sigma_{x_2}=\sigma_{x_3}\sigma_{x_2}\sigma_{x_3}$  is reduced.   
If $w=w(S)$ where $\displaystyle{S}$ is a (+)-admissible sequence of length $3$ on
$\displaystyle{(\Gamma, \Lambda)}$ for some $\displaystyle{\Lambda}$,
then either $S={x_2},{x_3},{x_2}$ or $S={x_3},{x_2},{x_3}$.  In the former case we must have $\displaystyle{a : x_1 \rightarrow x_2}$ in
$\gl$.  Then  in $(\Gamma,\displaystyle{\sigma_{x_3}\sigma_{x_2}\La})$ we have $\displaystyle{a : x_2 \rightarrow x_1}$ whence $\displaystyle{x_2}$ is not a sink, a contradiction.  If $S={x_3},{x_2},{x_3}$, the argument is the same, using $b$ instead of $a$.
\end{ex}

\begin{cor}\label{reduced} Suppose $\Gamma$ is not a Dynkin diagram of the type $A,D$, or $E$ and let $S\in\mathfrak S$.
\begin{itemize} 
\item[(a)] The word $w(S)  \in \mathcal W$ is reduced.
\item[(b)]  There exists a preprojective $k\gl$-module $M$ satisfying $S\sim S_M$.
\end{itemize}
\end{cor}
\begin{proof}  (a) By assumption, the finite-dimensional algebra $k\gl$ is of infinite representation type (see ~\cite{bgp}), whence there exist infinitely many  nonisomorphic indecomposable preprojective $k\gl$-modules $M$  ~\cite[VIII Proposition 1.16] {ars} and, by Theorem \ref{shrtstsq}(b), the corresponding sequences $S_M\in\mathfrak S$ are pairwise nonequivalent.  Since the poset $(\Gamma_0,\La)$ is finite, Proposition \ref{tails} implies that for a given $m\ge0$, there exists a sequence $S_M$ whose canonical form $T_1\dots T_q$ satisfies $m\le q$ and $T_i=K$  for $0<i\le m$ where $K$ is a complete sequence.  By Theorem \ref{reducedword}(a),  $w(S_M)$ is reduced whence so is $w(K^m)$ according to Remark \ref{SubseqReduced}.  For any $S\in\mathfrak S$, Proposition \ref{Prop:TFAE} implies $S\preccurlyeq K^r$ where $r$ is the size of $S$.  Using Remark \ref{SubseqReduced} again, we see that $w(S)$  is reduced.

(b) This is an immediate consequence of (a) and Theorem \ref{seqreducedwrd}.
\end{proof}

\begin{rmk}  In view of Theorem \ref{seqreducedwrd} and Corollary \ref{reduced}, for  a given $S\in\mathfrak S$ one may ask how to determine whether the word $w(S)$ is reduced;  and if yes, how to find a preprojective module $M$ satisfying $S\sim S_M$.  To handle these questions efficiently, one should write $S$ as the join of the smallest possible number of sequences $T_i\in\mathfrak P$ as explained in Proposition \ref{SeqJoinPrinc}; verify that each $w(T_i)$ is reduced using Theorem \ref{princseqreducedwrd}(c); and set $M$ to be the direct sum of $M_i$'s, where $M_i$ is the indecomposable preprojective $k\gl$-module obtained from $T_i$ according to Theorem \ref{shrtstsq}(d).
\end{rmk}

We now characterize  infinite Weyl groups in terms of reduced words.

\begin{thm}\label{infWeylGrp}  Let $A=(a_{ij})$ be an indecomposable symmetric generalized $n\times n$ Cartan matrix, and let $c=\sigma_{v_n}\dots\sigma_{v_1}$ be a Coxeter element of the Weyl group $\mathcal W$.   Then $\mathcal W$ is infinite if and only if for all $m\in\mathbb Z$, $\ell(c^m)=|m|n$.
\end{thm}
\begin{proof}  The sufficiency is clear.  For the necessity, note that there exists a unique  orientation $\La$ on $\Gamma$ for which  the quiver $\gl$ has no oriented cycles and $K=v_1,\dots,v_n$ is a (+)-admissible sequence on $\gl$ ~\cite[p. 8]{dr}.  Then  $c=w(K)$ and $c^m=w(K^m)$   for all $m\ge0$.  Since $\mathcal W$ is  infinite, ~\cite[Lemma 2.1, part 4), and Proposition 2.1]{bgp} say that $\Gamma$ is not a Dynkin diagram of the type A, D, or E.  By Corollary \ref{reduced}(a), $c^m$ is a reduced word.
\end{proof}

\section{(+)-admissible sequences and Coxeter-sortable elements}\label{c-sort}

The following definition quotes ~\cite[pp. 7-8]{rea}.
\begin{defn}\label{c-sort1} Fix an arbitrary Coxeter element $c=\sigma_{v_n}\dots\sigma_{v_1}$ in $\mathcal W$ and write a half-infinite sequence of vertices 
$$c^{\infty}={v_n},{v_{n-1}},\dots,{v_1},{v_n},{v_{n-1}},\dots,{v_1},{v_n},{v_{n-1}},\dots,{v_1},\dots$$  
The $c${\it -sorting word} for $w\in\mathcal W$ is the lexicographically first subsequence $v_{i_1},\dots,v_{i_s}$ of $c^{\infty}$ for which $\sigma_{v_{i_1}}\dots\sigma_{v_{i_s}}$  is a reduced word for $w$.  The $c$-sorting word can be interpreted as a sequence of subsets of $\Gamma_0$ by rewriting 
$$c^{\infty}={v_n},{v_{n-1}},\dots,{v_1}|{v_n},{v_{n-1}},\dots,{v_1}|{v_n},{v_{n-1}},\dots,{v_1}|\dots$$
where the symbol \lq\lq$|$" is called a {\it divider}.  The subsets in the sequence are the sets of vertices of the $c$-sorting word that occur between adjacent dividers.  This sequence contains a finite number of nonempty subsets, and if any subset is empty, then every later subset is also empty.  An element $w\in\mathcal W$ is $c${\it -sortable} if its $c$-sorting word defines a sequence of subsets that is decreasing under inclusion.
\end{defn}

\begin{prop}\label{c-sort} Let $K=v_1,\dots,v_n$ be a complete (+)-admissible sequence on $\gl$ and let $S\in\mathfrak S$. 
\begin{itemize}
\item[(a)] $c=\sigma_{v_n}\dots\sigma_{v_1}$ is a Coxeter element of $\mathcal W$.
\item[(b)]  If $S\sim S_M$ for some preprojective $k\gl$-module $M$, then  $w(S)^{-1}$ is a $c$-sortable element of $\mathcal W$. 
\end{itemize}
\end{prop}
\begin{proof} (a)  This is clear.

(b)  If $S=x_1,\dots,x_s$, then $S^t=x_s,\dots,x_1$  is a $(-)$-admissible sequence with respect to a suitable orientation, and $w(S)^{-1}=w(S^t)=\sigma_{x_s}\dots\sigma_{x_1}$.  By Theorem \ref{seqreducedwrd}, the word $w(S)$ is reduced, hence so is $w(S)^{-1}$.  By Proposition \ref{Prop:canonical}, $S\sim S_1S_2\dots S_r$ where each $S_i$ consists of distinct vertices and $\mathrm{Supp}\, S_{i+1}\subset \mathrm{Supp}\, S_{i}$.  Then $w(S^t)=w(S_r^t)\dots w(S_1^t)$ where $\mathrm{Supp}\, S_{i+1}^t\supset \mathrm{Supp}\, S_{i}^t$. 
\end{proof}

\end{document}